# Interactions of Andronov–Hopf and Bogdanov–Takens Bifurcations

William F. Langford and Kaijun Zhan

Abstract. A codimension-three bifurcation, characterized by a pair of purely imaginary eigenvalues and a nonsemisimple double zero eigenvalue, arises in the study of a pair of weakly coupled nonlinear oscillators with $\mathbb{Z}_2 \oplus \mathbb{Z}_2$ symmetry. The methodology is based on Arnold's ideas of versal deformations of matrices for the linear analysis, and Poincaré normal forms for the nonlinear analysis of the system. The stratified subvariety of primary bifurcations of codimensions one and two is identified in the parameter space. The analysis reveals different types of solutions in the state space, including equilibria, limit cycles, invariant tori and the possibility of homoclinic chaos. A mechanism is identified for energy transfer without strong resonance between two oscillation modes with widely separated frequencies.

## 1. Introduction

This work is motivated by the study of weakly-coupled nonlinear oscillators which may be modeled by second-order differential equations of the general form

$$\begin{aligned} \ddot{x} + \omega_1^2 x + f_1(x, \dot{x}, y, \dot{y}, \mu) = 0 \\ \ddot{y} + \omega_2^2 y + f_2(x, \dot{x}, y, \dot{y}, \mu) = 0. \end{aligned} \tag{1.1}$$

Here $f_1$ and $f_2$ are holomorphic functions of the state variables $x, \dot{x}, y, \dot{y}$ and parameters $\mu \in \mathbb{R}^k$. The functions are higher order in that they contain only quadratic and higher order terms in these state variables and parameters, which are assumed to be small, and $f_j(0,0,0,0,\mu) \equiv 0,\ j = 1,2$. The parameters $\mu$ represent physically relevant quantities such as damping, detuning and the coupling between the oscillators; the role of these parameters will be clarified further below. In (1.1), the natural frequencies of the two oscillators, when linearized at $(x, \dot{x}, y, \dot{y}) = (0,0,0,0)$, are given by $\omega_{1,2} \geq 0$. Systems of the form (1.1) occur frequently in the study of mechanical systems, electrical circuits, the biological sciences and elsewhere.

A major concern in the study of coupled oscillators is the understanding and control of *resonance* phenomena. V.I. Arnold has made influential contributions to the general theory of resonances, for example [**5, 8, 9**]. The ratio of the natural frequencies plays a leading role in determining resonances in systems (1.1) of two oscillators. The cases $\frac{\omega_1}{\omega_2} = \frac{p}{q}$, where $p$ and $q$ are small positive integers, $p +$

1991 *Mathematics Subject Classification.* Primary 58F14; secondary 34C20.

Research supported by the Natural Sciences and Engineering Research Council, Canada, The Fields Institute for Research in Mathematical Sciences, Toronto ON and the Institute for Mathematics and its Applications, Minneapolis MN.

$q \leq 5$, are called *strong resonances*. The remaining rational cases are called *weak resonances*, and cases with the ratio $\frac{\omega_1}{\omega_2}$ irrational are called *nonresonant*. As the names suggest, coupled oscillators interact strongly at or near strong resonances, and interact weakly, if at all, in the other cases. This is explained in part by the fact that in the Poincaré normal form corresponding to these coupled oscillator equations, terms of degree $p+q$ in the state variables remain. These are called resonant terms. Locally, these terms have a strong influence when $p+q$ is small and their influence decreases as $p+q$ increases.

This paper addresses the case of coupled oscillators (1.1) with natural frequencies which are widely separated; we assume $\omega_1$ very small and $\omega_2$ very large, relative to each other. On the basis of the above classification, this would not be considered a case of strong resonance. To fix magnitudes, rescale time in (1.1) so that $\omega_2 = 1$. Then our assumption is $\omega_1 \to 0$. Equivalently, we consider natural frequencies at or near the ratio $\frac{\omega_1}{\omega_2} = \frac{0}{1}$. Therefore, we refer to this as the case of $0:1$ resonance.

In this paper, system (1.1) is viewed as a perturbation of a singular system with a bifurcation of codimension three, arising through the coalescence of two simpler bifurcations. The low frequency mode, in the limit $\omega_1 \to 0$, is a Bogdanov–Takens bifurcation [**11, 41**]. This bifurcation has codimension two. The high frequency mode is close to a nondegenerate Andronov–Hopf bifurcation, which has codimension one. Historically, this bifurcation has its origins in the work of Poincaré [**37**, sections 51–52,316] more than a century ago, and was studied extensively by Andronov and coworkers [**1, 2, 3, 4**] starting in the 1920's. Arnold [**9**, section 33.A,B] presents a pellucid account of the Poincaré–Andronov theory. An excellent account of the approach of Hopf and later developments is given in [**30**]. Degenerate cases of Hopf's theorem are explored in [**17**]. Although many different names have been used for this theory over the years, in this paper we choose the name Andronov–Hopf bifurcation, which is suggestive of its rich history.

Using Arnold's theory of versal deformations of matrices [**6**], we verify that the singular system combining these two classical bifurcations has linear part with codimension three. Then we investigate via Poincaré normal forms the implications of this codimension-three singularity for nonlinear systems. The analysis reveals, in addition to features associated with the original two classical bifurcations, new phenomena arising from the nonlinear mode-interactions.

Cases of bifurcations involving a *simple* zero eigenvalue with a pair of purely imaginary eigenvalues (and no other eigenvalues with zero real part) have been studied by Arnold [**9**], Langford [**24, 25**], Guckenheimer and Holmes [**18**], Khazin and Shnol [**23**], Chow et al. [**14**], Scheurle and Marsden [**39**] and others cited in these references. Less is known about the case of a nonsemisimple *double* zero eigenvalue together with an imaginary pair. Some early analytical results and applications have been presented by Moson [**31, 32**], Yu and Huseyin [**45**], Georgiou et al. [**15**], Nagata and Namachchivaya [**33**], and Nayfeh et al. [**34, 35, 36**].

In a series of papers [**34, 35, 36**], Nayfeh and coworkers have investigated a system of two weakly coupled oscillators with widely separated frequencies similar to (1.1). They point out an interesting paradox. In both mechanical experiments and a mathematical model, when the system is excited near a high natural frequency, large low-frequency responses accompanied by slow modulations of the amplitude and phase of the high-frequency mode are observed. Thus, they found that there is a strong interaction between these two modes, and a transfer of energy from high to low frequency modes, even though by the above classification it is not a case of strong resonance. Their work provided additional motivation for the present paper.

To fix ideas, throughout this paper it is assumed that each of the oscillators in (1.1) has an odd reflectional symmetry about the rest position. We extend this symmetry to the coupled system by assuming that the system commutes with the $\mathbb{Z}_2 \oplus \mathbb{Z}_2$ symmetry group generated by two reflection operators

$$\begin{aligned} f_1(-x,-\dot{x},y,\dot{y},\mu) &= -f_1(x,\dot{x},y,\dot{y},\mu) \\ f_2(-x,-\dot{x},y,\dot{y},\mu) &= f_2(x,\dot{x},y,\dot{y},\mu) \\ f_1(x,\dot{x},-y,-\dot{y},\mu) &= f_1(x,\dot{x},y,\dot{y},\mu) \\ f_2(x,\dot{x},-y,-\dot{y},\mu) &= -f_2(x,\dot{x},y,\dot{y},\mu)\ . \end{aligned} \tag{1.2}$$

This symmetry is very common in applications; for example, the classical pendulum equation, Duffing's equation and Van der Pol's equation all have odd symmetry, as do the examples in [**33, 34, 35, 36**]. The symmetry (1.2) forbids linear coupling terms (i.e. linear $x, \dot{x}$ terms in the $\ddot{y}$-equation or vice-versa) but permits nonlinear coupling (such as $x^2y$ in the $\ddot{y}$-equation). Systems with less symmetry that (1.2) may also be studied by the methods of this paper.

Before beginning to analyze (1.1), we rewrite it with the linear part in a more explicit form. Equation (1.3) displays the most general linear terms in the state variables, which are consistent with the symmetries (1.2) and the time-rescaling. (Recall that $t$ has been rescaled so that $\omega_2 = 1$). The functions $\hat{f}_j$, which replace the former functions $f_j$, consist of all the remaining nonlinear terms in the state variables $(x,\dot{x},y,\dot{y})$, with coefficients which are functions of the parameters $\mu$. Since the functions $\hat{f}_j$ are assumed to commute with the two reflection operators (1.2), they are cubic or higher order in the state variables.

$$\begin{aligned} \ddot{x} + \delta_1\dot{x} + \Omega x + \hat{f}_1(x,\dot{x},y,\dot{y},\mu) &= 0 \\ \ddot{y} + \delta_2\dot{y} + y + \hat{f}_2(x,\dot{x},y,\dot{y},\mu) &= 0 \end{aligned} \tag{1.3}$$

In the linear terms of (1.3), $\delta_1, \delta_2$ are linear damping parameters and $\Omega = (\omega_1/\omega_2)^2$ is a small tuning parameter which from here on is allowed to be negative or zero; therefore, $\Omega$ is better thought of as a restoring force than a frequency ratio. In the engineering literature, a change in sign of $\Omega$ from positive to negative is associated with a transition from "flutter" to "divergence". These physical parameters are all assumed to be small. Recognizing the physical importance of these three parameters and choosing to ignore the possible effects of other parameters, we formally define the parameter vector $\mu$ to have precisely these three components

$$\mu \equiv (\Omega, \delta_1, \delta_2)\ . \tag{1.4}$$

The appropriateness of this choice will become clearer after the presentation of Arnold's versal deformation, in the next section.

With the definition

$$\mathbf{y} \equiv \begin{pmatrix} y_1 \\ y_2 \\ y_3 \\ y_4 \end{pmatrix} \equiv \begin{pmatrix} x \\ \dot{x} \\ y \\ \dot{y} \end{pmatrix}$$

system (1.3) can be written as the following system of first order equations, with linear part which is a perturbation of a matrix $\mathbf{A}_0$ in real Jordan form.

$$\dot{\mathbf{y}} = [\mathbf{A}_0 + \mathbf{C}(\mu)]\,\mathbf{y} + \hat{\mathbf{f}}\ , \tag{1.5}$$

$$\mathbf{A}_0 = \begin{pmatrix} 0 & 1 & 0 & 0 \\ 0 & 0 & 0 & 0 \\ 0 & 0 & 0 & 1 \\ 0 & 0 & -1 & 0 \end{pmatrix}, \mathbf{C}(\mu) = \begin{pmatrix} 0 & 0 & 0 & 0 \\ -\Omega & -\delta_1 & 0 & 0 \\ 0 & 0 & 0 & 0 \\ 0 & 0 & 0 & -\delta_2 \end{pmatrix}, \hat{\mathbf{f}} = \begin{pmatrix} 0 \\ -\hat{f}_1 \\ 0 \\ -\hat{f}_2 \end{pmatrix}.$$

Here $\mathbf{C}(\mu)$ is the linear perturbation matrix, depending only on the three small physical parameters $\mu = (\Omega, \delta_1, \delta_2)$.

We may simplify notation by introducing a mixed set of real/complex state variables with $(y_1, y_2) \in \mathbb{R}^2$ defined as above, and $z \in \mathbb{C}$ defined by

$$z = y_3 - iy_4, \ \bar{z} = y_3 + iy_4 \ . \tag{1.6}$$

(The reason for the nonstandard choice of signs of the imaginary parts of $z, \bar{z}$ is the following. In the phase plane of a second order harmonic oscillator equation, the coordinate axes are normally chosen to be $(x, y) = (x, \dot{x})$, which forces oscillations about the origin to have a clockwise orientation. In the complex plane, the convention is that a counterclockwise rotation is positive, as for the differential equation $\dot{z} = iz$. The transformation (1.6) effects this reversal of orientation, between the real and complex phase planes.) Then system (1.5) is transformed to the system

$$\begin{aligned} \dot{y}_1 &= y_2 \\ \dot{y}_2 &= -\Omega y_1 - \delta_1 y_2 - \hat{f}_1 \\ \dot{z} &= (i - \tfrac{1}{2}\delta_2)z + \tfrac{1}{2}\delta_2 \bar{z} + i\hat{f}_2 \end{aligned} \tag{1.7}$$

together with a fourth equation for $\dot{\bar{z}}$ obtained by conjugation of the equation for $\dot{z}$ given here.

The matrix $\mathbf{A}_0$ in (1.5) is recognized as being in real Jordan form, with the four eigenvalues $0^2$ (nonsemisimple) and $\pm i$ (imaginary pair). The corresponding differential equations (1.5) for small amplitudes represent a perturbation of a system of two independent oscillators, one with coordinates $(x, \dot{x})$ and natural frequency $\omega_1 \approx 0$, and the other with coordinates $(y, \dot{y})$ and natural frequency $\omega_2 = 1$. However, the results of this paper apply more generally than this. Consider any parametrized two degree-of-freedom system, which may have coefficient matrix far from $\mathbf{A}_0$; but suppose that, for some value of the parameters, the linearization at an equilibrium has eigenvalues $\{0^2, \pm i\}$ (after rescaling time). Then the linear part can be transformed to $\mathbf{A}_0$ by an invertible linear change of coordinates (similarity transformation), which may be far from the identity. In this case, the state variables in (1.5) are the *normal mode* coordinates of the original system, and the symmetries (1.2) are symmetries of the normal modes. The state variables in this paper should be understood in that light.

Even more generally, the results of this paper apply to systems of ordinary differential equations of dimension higher than four, and to classes of infinite dimensional partial differential equations or functional differential equations, satisfying appropriate spectral conditions. If the linearization at an equilibrium solution satisfies standard compactness and hyperbolicity assumptions on the linear operator and its spectrum, then the Center Manifold Theorem applies. In this case, if the spectrum contains the eigenvalues $\{0^2, \pm i\}$, and the complementary spectrum is bounded away from the imaginary axis in the left half-plane, then the study of long-time dynamics of the full system reduces to the four dimensional system considered here. Such higher dimensional generalizations are not pursued further.

The paper is organized as follows. In the next section, a versal deformation of the linear part of the singular vector field is obtained, using ideas of Arnold, and the codimension is computed to be three. Section 3 presents the Poincaré normal form for this system. The original four-dimensional system is reduced to a three-dimensional one, exploiting the $S^1$ normal form symmetry arising from the Andronov–Hopf bifurcation. In Section 4, the generic bifurcations of codimension one and two are explored. A rich variety of solutions, including stationary, periodic and quasiperiodic (tori) are found. The main results of this paper include the stratified subvariety of primary bifurcations, and formulas for primary and secondary bifurcations of equilibria, periodic solutions and tori. In Section 5, the relevance of these results to the energy-transfer phenomenon of Nayfeh [**34, 35, 36**] is considered by means of an example, and directions of further research are indicated in Section 6.

## 2. Linear Analysis

**Definition.** Let $\mathbb{C}^{n^2}$ denote the set of $n \times n$ complex matrices. A *family $A$ of matrices on the base $\mathcal{M}$* is a holomorphic mapping $A : \mathcal{M} \to \mathbb{C}^{n^2}$, where $\mathcal{M}$ is a neighbourhood of the origin in the parameter space $\mathbb{C}^k$. For brevity, we call this the *family $A$*. The germ of a family $A$ at the point $0 \in \mathcal{M}$ is called a *deformation* of the matrix $A(0)$. The elements $\mu \in \mathcal{M}$ are called *parameters*.

An example of a family of real matrices is given by $\mathbf{A}_0 + \mathbf{C}(\mu)$ in equation (1.5) and a complex family is given by the linear part of (1.7). The algebra is simpler if we work with complex matrices. The real Jordan matrix $\mathbf{A}_0$ of (1.5) has the complex Jordan form

$$\mathbf{A}_c = \begin{pmatrix} 0 & 1 & 0 & 0 \\ 0 & 0 & 0 & 0 \\ 0 & 0 & i & 0 \\ 0 & 0 & 0 & -i \end{pmatrix} . \tag{2.1}$$

Two main questions are addressed in this section. The first question is how to find a canonical matrix, which is as simple as possible, and to which the linear part of all coupled-oscillator systems as described in Section 1 can be reduced by a linear change of coordinates. An obvious candidate is the Jordan form; however, this is not a good choice because the reduction of a matrix to its Jordan form is not a stable operation when there are multiple eigenvalues, as is the case here. Although every member of a family containing $\mathbf{A}_c$ may be reduced to a Jordan form, the Jordan form does not in general depend continuously on the parameters, since there are matrices arbitrarily close to $\mathbf{A}_c$ for which the Jordan form is diagonal. The mapping which transforms matrices to Jordan form is discontinuous. Therefore, we seek a different normal form, to which members of any family containing $\mathbf{A}_c$ may be transformed smoothly.

The second question relates to the choice of parameters $\mu \in \mathcal{M}$. What is the minimum set of parameters needed, and how can we know if any essential parameters are missing? In particular, are there systems with linear parts close to $\mathbf{A}_0$ which are not attainable by variations of the physical parameters in the given perturbation matrix $\mathbf{C}(\mu)$ in (1.5)? If we have in fact overlooked some essential parameters, then any study of the system in Section 1 is likely to miss important phenomena. We would like to have assurance that the parameter set in our family of matrices is sufficiently rich that *all* perturbations of $\mathbf{A}_0$ are included, and that there are no redundant parameters. On the other hand, it may happen that due to the mechanical structure of equations (1.3), the parameters and forces can enter only in certain ways and generic perturbations are unrealistic for the system. In that case, a comparison with the generic case is useful, to uncover any previously hidden constraints in our model.

Both of these classes of questions were answered by Arnold [**6, 7, 9**], who introduced a generalization of the Jordan form of a matrix $A$, with the minimum number of essential parameters, and to which all members of any family of $A$ can be transformed smoothly. This generalization of the Jordan normal form may be called the *Arnold normal form* of a matrix. The remainder of this section summarizes Arnold's construction of this normal form and shows how it answers the two questions above, for a pair of coupled oscillators with linear part as in (1.5). A general exposition of these ideas of Arnold is given in [**43**, pages 305–320].

The reduction to Jordan form of any matrix $A$ is accomplished by a similarity transformation. Consider the Lie group $\mathrm{GL}(n,\mathbb{C})$ of all nonsingular matrices $S \in \mathbb{C}^{n^2}$. Then the group of similarity transformations on $\mathbb{C}^{n^2}$ has the elements $\mathrm{Ad}_S,\ S \in \mathrm{GL}(n,\mathbb{C})$, defined by

$$\mathrm{Ad}_S(A) \equiv SAS^{-1}, \quad A \in \mathbb{C}^{n^2},\ S \in \mathrm{GL}(n,\mathbb{C}) \ . \tag{2.2}$$

The *group orbit* of a given $A$ is

$$\mathcal{O}(A) = \left\{ \mathrm{Ad}_S(A) = SAS^{-1} \ \mid \ S \in \mathrm{GL}(n,\mathbb{C}) \right\} \tag{2.3}$$

and this is a smooth submanifold of $\mathbb{C}^{n^2}$. The orbit of $A$ thus consists of all matrices similar to $A$. It is well known that similar matrices share the same eigenvalues, with the same multiplicities and the same Jordan blocks (perhaps permuted). A convenient orbit representative for any given $A$ is its Jordan form, with Jordan blocks ordered in decreasing size, for example $\mathbf{A}_c$ in (2.1).

Now, combining the concepts of deformations and similarity, we say that two deformations $B(\mu)$ and $A(\mu)$ of the same matrix $B(0) = A(0) = A_0$ are *equivalent* if there exists a deformation $C(\mu)$ of the identity matrix, for $\mu \in \mathcal{M}$, such that

$$B(\mu) = C(\mu)A(\mu)C(\mu)^{-1}. \tag{2.4}$$

Next, this equivalence is extended to allow two different parametrizations for $B$ and $A$. Let $\mathcal{M}$ and $\mathcal{N}$ be neighbourhoods of the origin in parameter spaces $\mathbb{C}^k$ and $\mathbb{C}^l$, respectively. Let $\varphi : \mathcal{M} \to \mathcal{N}$ be holomorphic at 0, and $\varphi(0) = 0$. Let $A(\alpha)$ be a family of matrices on $\mathcal{N}$. Then the family (deformation) *induced from $A$ by $\varphi$* is the family (deformation) $\varphi \circ A$

$$(\varphi \circ A)(\mu) = A(\varphi(\mu)).$$

**Definition.** A deformation $A(\alpha)$ of $A_0$ is called *versal* if every deformation $B(\mu)$ of $A_0$ is equivalent to a deformation induced from $A(\alpha)$:

$$B(\mu) = C(\mu)A(\varphi(\mu))C^{-1}(\mu), \quad C(0) = I,\ \varphi(0) = 0. \tag{2.5}$$

A versal deformation is *universal* if mapping $\varphi$ is determined uniquely by deformation $B(\mu)$. It is *miniversal* if the dimension of the parameter space $\mathcal{N}$ is minimal among all versal deformations. This minimum dimension is called the *codimension* of the matrix $A_0$. The parameters in a miniversal deformation are also called *unfolding parameters*, and a (uni-, mini-)versal deformation is sometimes called a *(uni-, mini-)versal unfolding.*

It is clear that such a versal deformation exists and an upper bound on the codimension is $n^2$, since one could take as parameters all the elements of an $n \times n$ matrix. A lower bound is given by $n$, for the case of distinct eigenvalues and a diagonal Jordan form. In general, we need a manifold transverse to the group orbits. To this end, compute tangents to the group orbits in (2.3), for fixed $A = A_0$. The derivative of $\mathrm{Ad}_S(A_0) = SA_0S^{-1}$, with respect to $S$ at the identity $I$, for fixed $A_0$, is easily found to be the linear operator defined by commutation with $A_0$

$$DAd_I : \mathbb{C}^{n^2} \to \mathbb{C}^{n^2}, \quad DAd_I S = [S, A_0] \equiv SA_0 - A_0 S. \tag{2.6}$$

The *centralizer* $Z_A$ of matrix $A \in \mathbb{C}^{n^2}$ is the set of all matrices which commute with $A$

$$\mathrm{Z}_A = \{ \ B \ | \ [B, A] \equiv BA - AB = 0\} .$$

Thus, the nullspace of $D\mathrm{Ad}_I$ is the centralizer of $A_0$. The rank of $D\mathrm{Ad}_I$ is the dimension of the group orbit of $A_0$. Since $D\mathrm{Ad}_I$ is square, the dimension of the centralizer is the codimension of the group orbit, which is the minimum number of parameters in a versal deformation. This codimension was given by Arnold [**6**] in the following formula

$$d = \sum_{i=1}^{s} \sum_{j=1}^{l_i} (2j-1) n_{ij}(\lambda_i) \tag{2.7}$$

where $\lambda_i, \ i = 1, \ldots, s$ are distinct eigenvalues of $A_0$, $l_i$ is the number of Jordan blocks of $\lambda_i$ and $n_{i1}(\lambda_i) \geq \cdots \geq n_{il_i}(\alpha_i)$ are the corresponding orders of the Jordan blocks for each $\lambda_i$.

Furthermore, with respect to the Hermitian scalar product $\langle A, B \rangle = \mathrm{tr}(AB^*)$, the orthogonal complement of the tangent to the orbit of $A_0$ in the tangent space $\mathbb{C}^{n^2}$ is given by $\{B^* \mid B \in \mathrm{Z}_{A_0}\}$. Here $B^*$ denotes the conjugate transpose of $B$. This is the form that the Fredholm Alternative takes in this setting. From this, Arnold obtains explicitly the versal deformations, see [**6, 7, 9, 43**].

Applying Arnold's method to the Jordan matrix $\mathbf{A}_c$ above gives the centralizer $\mathrm{Z}_{\mathbf{A}_c} = \{C(\alpha)\}$ and versal deformation $\mathbf{A}_c(\alpha)$ explicitly as the families:

$$C(\alpha) = \begin{pmatrix} \alpha_2 & \alpha_1 & 0 & 0 \\ 0 & \alpha_2 & 0 & 0 \\ 0 & 0 & \alpha_3 & 0 \\ 0 & 0 & 0 & \alpha_4 \end{pmatrix}, \tag{2.8}$$

$$\mathbf{A}_c(\alpha) = \begin{pmatrix} \alpha_2 & 1 & 0 & 0 \\ \alpha_1 & \alpha_2 & 0 & 0 \\ 0 & 0 & i + \alpha_3 & 0 \\ 0 & 0 & 0 & -i + \alpha_4 \end{pmatrix}. \tag{2.9}$$

Here $\alpha_j, \ j = 1, \ldots, 4$ are arbitrary complex parameters. This versal deformation is not unique. Rather than the above orthogonal complement to the tangent space, we could choose any transversal subspace of minimum dimension, and obtain a miniversal deformation. One convenient choice replaces the $\alpha_2$ entry in the upper left corner of $\mathbf{A}_c(\alpha)$ by 0, and thus reduces the number of nonzero entries by one. We adopt this choice.

Furthermore, since the original system is real, it is desirable to decomplexify $\mathbf{A}_c(\alpha)$. The result of decomplexification (for this choice) is the real family

$$\mathbf{A}_0(\alpha) = \begin{pmatrix} 0 & 1 & 0 & 0 \\ \alpha_1 & \alpha_2 & 0 & 0 \\ 0 & 0 & \alpha_3 & -1-\alpha_4 \\ 0 & 0 & 1+\alpha_4 & \alpha_3 \end{pmatrix} \tag{2.10}$$

where $\alpha_j,\ j = 1, \ldots, 4$ have been redefined as real parameters. We call the family in equation (2.10) the real Arnold normal form of the original matrix $\mathbf{A}_0$. It has the desired versality property, that any deformation $B(\mu)$ of $\mathbf{A}_0$ is equivalent to a deformation induced from the family $\mathbf{A}_\alpha$. In other words, any such $B(\mu)$ can be obtained by means of holomorphic mappings of parameters combined with similarity transformations, as given in (2.5). Among all such versal deformations, it has the minimum number of parameters (i.e. is miniversal).

Finally, recall that our matrix is the linear part of a differential equation, in which a rescaling of time has the effect of multiplying the matrix by a positive real number. We choose to rescale time to keep the imaginary part of the complex eigenvalue equal to $\pm i$. This has the effect of making $\alpha_4 = 0$ in (2.10) and rescaling the other parameters. Therefore, in this context the codimension is 3, and the versal deformation is

$$\mathbf{A}_0(\alpha) = \begin{pmatrix} 0 & 1 & 0 & 0 \\ \alpha_1 & \alpha_2 & 0 & 0 \\ 0 & 0 & \alpha_3 & -1 \\ 0 & 0 & 1 & \alpha_3 \end{pmatrix} \tag{2.11}$$

where now $\alpha_j,\ j = 1, 2, 3$ are redefined real unfolding parameters.

It remains to establish a connection between the unfolding parameters $\alpha_j$ in Arnold's versal deformation (or normal form) (2.11) and the physical parameters in the original system (1.5). Since the physical parameters can enter only in certain ways, determined by the mechanical structure of the system, it is not clear that the notion of versality for general differential equations is appropriate for this class. If in fact the physical parameters can not be transformed to the Arnold versal deformation parameters, then two possibilities must be considered: either there is a constraint on the physical system which should be made explicit, or an allowable parameter has been overlooked.

For the present example, it is easy to find an explicit relationship between the physical and unfolding parameters. This may be done by comparing coefficients of the characteristic polynomials. The relationship is found to be

$$\begin{aligned} \alpha_1 &= -\Omega + \cdots \\ \alpha_2 &= -\delta_1 + \cdots \\ \alpha_3 &= -\tfrac{1}{2}\delta_2 + \cdots \end{aligned} \tag{2.12}$$

where $\cdots$ represents higher order terms. Since this mapping is invertible at the origin, we have proven the following.

**Proposition.** The linear part of the physical model in (1.5), with parameters $\{\Omega, \delta_1, \delta_2\}$, is a miniversal deformation in the sense of Arnold. The three physical parameters are related to the unfolding parameters in (2.11) by equations (2.12).

This justifies the choice of three physical parameters $(\Omega, \delta_1, \delta_2)$ in equation (1.3), in the following sense. Although there may be other physical parameters which affect the behaviour of the model, those other parameters can not lead to any new local behaviour which is not already included in the universal deformation (2.11). Thus, the 3-parameter family of matrices (2.11) is the appropriate starting point for the study of phenomena associated with coupled oscillators near 0 : 1 resonance. As observed by Arnold in [**9**, page 221]:

> ... the natural object of study is not the individual object (say, a vector field with a complicated singular point), but a family so large that the singularity of the type under consideration does not disappear under a small deformation of the family. This simple argument of Poincaré shows the futility of such a large number of studies in the theory of differential equations and in other areas of analysis that it is always somewhat dangerous to mention it.

### 3. Poincaré Normal Form

The Poincaré normal form may be thought of as an extension of the ideas of the previous section from the linear case to the case of *nonlinear* vector fields. Arnold gives an excellent account of the theory and computation of the Poincaré normal form in [**9**, Chapter 5]. Other useful references are [**10, 12, 14, 18, 20, 22, 23, 28, 43**].

The basic idea of Poincaré is to perform a sequence of near-identity transformations of the form

$$\mathbf{x} = \mathbf{y} + \mathbf{h}_k(\mathbf{y}), \quad k = 2, 3, \dots \tag{3.1}$$

where $\mathbf{h}_k$ is a homogeneous polynomial of degree $k$. For each $k$, $\mathbf{h}_k$ is chosen to eliminate as many as possible of the terms of degree $k$ in the vector field. Terms which remain lie in a complement to the tangent space, at the identity, of the manifold of group orbits corresponding to (3.1), just as in the linear case of Section 2. The Poincaré normal form is not unique, since it depends on this choice of complement. Furthermore, the resulting power series does not converge in general; nor does the sequence of near-identity transformations (3.1) converge. This lack of convergence is often not a problem, since in many cases the phenomena of interest are determined by a finite truncation of the Poincaré normal form, of low degree.

For the system of differential equations in this paper as given by (1.7) with $\mu = 0$ , the Poincaré normal form was first computed by Moson [**31, 32**]; see also [**10, 12**]. The Poincaré normal form inherits the symmetry property (1.2), and can be chosen to have an additional *normal form symmetry* which is the $\mathrm{S}^1$ symmetry inherent in the Andronov–Hopf bifurcation [**20**]. The resulting normal form may be written (formally) as

$$\begin{pmatrix} \dot{y}_1 \\ \dot{y}_2 \\ \dot{z} \end{pmatrix} = \begin{pmatrix} 0 & 1 & 0 \\ 0 & 0 & 0 \\ 0 & 0 & i \end{pmatrix} \begin{pmatrix} y_1 \\ y_2 \\ z \end{pmatrix} + \begin{pmatrix} \phi_2 y_1 \\ \phi_1 y_1 + \phi_2 y_2 \\ [\phi_3 + i\phi_4]\, z \end{pmatrix} \tag{3.2}$$

to which may be added the conjugate of the last equation, for $\dot{\bar{z}}$. Here $\phi_j,\ j = 1, \dots, 4$, are real valued formal power series (which may not converge) in the Hilbert basis of invariant polynomials $\{\ y_1^2,\ \ z\bar{z} = y_3^2 + y_4^2\ \}$, i.e.

$$\phi_j = \phi_j(y_1^2, z\bar{z}), \ \ j = 1, \dots, 4.$$

Notice the formal resemblance between (3.2) and the complex Arnold normal form in the linear case, given by (2.9). In fact, if we were to allow the terms involving $\phi_j$ in the vector on the right-hand side of (3.2) to have linear parts, then the linear part of the Poincaré normal form would be the Arnold normal form of its linearization.

The Poincaré normal form (3.2) can be further simplified by the near-identity transformation:

$$(y_1, y_2, z) \to (y_1, y_2 + \phi_2 y_1, z),$$

which eliminates the nonlinear terms in the $\dot{y}_1$ equation of (3.2). Now, the Poincaré normal form matches the Arnold normal form as given in (2.10), in the sense that we can identify

$$\alpha_j \equiv \phi_j(0,0), \quad j = 1, \ldots, 4 \ .$$

Combining in one equation the versal deformation of the linearization, i.e. the Arnold normal form given by (2.11), and this simplified Poincaré normal form, yields

$$\begin{pmatrix} \dot{y}_1 \\ \dot{y}_2 \\ \dot{z} \end{pmatrix} = \begin{pmatrix} 0 & 1 & 0 \\ \alpha_1 & \alpha_2 & 0 \\ 0 & 0 & i+\alpha_3 \end{pmatrix} \begin{pmatrix} y_1 \\ y_2 \\ z \end{pmatrix} + \begin{pmatrix} 0 \\ \phi_1 y_1 + \phi_2 y_2 \\ [\phi_3 + i\phi_4]\, z \end{pmatrix} \tag{3.3}$$

where $y_1, y_2 \in \mathbb{R}$, $z \in \mathbb{C}$, and $\phi_j = \phi_j(y_1^2, z\bar{z})$, $\phi_j(0,0) = 0$, $j = 1, \ldots, 4$, and we have rescaled time to make $i\alpha_4 = 0$.

Finally, we decomplexify the $z$-equation to real form, but in polar rather than cartesian coordinates, to better exploit the normal form symmetry. Let

$$z \equiv re^{i\theta} \ .$$

Then (3.3) can be written

$$\begin{aligned} \dot{y}_1 &= y_2 \\ \dot{y}_2 &= [\alpha_1 + \phi_1]y_1 + [\alpha_2 + \phi_2]y_2 \\ \dot{r} &= [\alpha_3 + \phi_3]r \\ \dot{\theta} &= 1 + \phi_4 \ . \end{aligned} \tag{3.4}$$

Note that $\theta$ does not appear in the right hand sides of equations (3.4), so the $\dot{\theta}$-equation is decoupled from the other three. This is due to the $S^1$ symmetry of the normal form, induced by the Andronov–Hopf bifurcation in $z$, and it effectively reduces the problem from four dimensions to three. Given a solution $(y_1(t), y_2(t), r(t))$ of the first three equations of (3.4), then the solution $\theta(t)$ is given by

$$\theta(t) = \theta_0 + t + \int_0^t \phi_4(y_1^2(s), r^2(s))\, ds \ .$$

All but exotic behaviour of this system is already determined by a cubic truncation of the three-dimensional Poincaré normal form. This can be written

$$\begin{aligned} \dot{y}_1 &= y_2 \\ \dot{y}_2 &= [\alpha_1 + \phi_{1,10}y_1^2 + \phi_{1,01}r^2]y_1 + [\alpha_2 + \phi_{2,10}y_1^2 + \phi_{2,01}r^2]y_2 \\ \dot{r} &= [\alpha_3 + \phi_{3,10}y_1^2 + \phi_{3,01}r^2]r \end{aligned} \tag{3.5}$$

where $\phi_{j,kl}$ is the coefficient of the term $y_1^{2k}r^{2l}$, in the power series expansion of $\phi_j(y_1^2, r^2)$. The $\theta$ equation has been omitted in (3.5), but has a similar expansion. The system (3.5) is explored in the next section.

By standard arguments, hyperbolic solutions of (3.5)are preserved as solutions of (3.4), and in fact correspond to higher dimensional solutions of the original four dimensional system. For example, hyperbolic equilibria of (3.5) correspond to hyperbolic periodic orbits of the original system. However, as has been shown in [**18, 25, 26, 43**], in special cases the higher order remainder or "tail" of the normal form can dramatically change the dynamics from the behaviour of a truncated normal form such as (3.5), since the tail breaks the normal form symmetry which led to the decoupling of the $\dot{\theta}$ equation. For example, when (3.5) has homoclinic or heteroclinic orbits, then the tail effects can produce homoclinic tangles, Smale horseshoes and chaotic dynamics. The singular nature of these normal form perturbations means that splitting of separatrices will be associated with exponentially small splitting problems, see [**38, 19, 40**]. These considerations will be left for investigation elsewhere.

## 4. Bifurcations of Codimension One and Two

Identifying the collection of all possible phase portraits of this system near the codimension three bifurcation point, in the sense of all equivalence classes with respect to topological equivalence, remains an open problem. However, the *stratified subvariety* (as defined in [**9**, page 227]) of bifurcations from the trivial solution in the parameter space $\mathbb{R}^3$ takes a remarkably simple form, and is presented in this section.

The stability of the trivial equilibrium solution is determined by the eigenvalues of the linear part $\mathbf{A}_0(\alpha)$, which are found to be

$$\lambda_{1,2} = \frac{\alpha_2}{2} \pm \sqrt{\alpha_1 + \left(\frac{\alpha_2}{2}\right)^2}, \quad \lambda_{3,4} = \alpha_3 \pm i \ . \tag{4.1}$$

All four eigenvalues have negative real part if and only if $\alpha_j < 0$, $j = 1,2,3$. Therefore, the trivial equilibrium solution is asymptotically stable, strictly inside the negative octant of the parameter space $\mathcal{N} \subset \mathbb{R}^3$, see Figure 1. As the point $(\alpha_1, \alpha_2, \alpha_3)$ leaves the negative octant, bifurcations occur. These bifurcations from the trivial equilibrium solution are called *primary bifurcations*. There are seven possibilities, see Figure 1: a point may leave by crossing through one of the three coordinate planes (codimension one), through one of the three coordinate axes (codimension two) or through the origin (codimension three). A generic exit is through a plane, codimension one. However, interesting and important mode-interactions occur around the intersections of the planes. These are best understood by studying the codimension two bifurcations on the coordinate axes. All six of these codimension one and two bifurcations have been studied previously, and are well understood in general. Their manifestations in the present system are described in the following subsections.

**4.1. Primary Bifurcations of Codimension One.** The three primary bifurcations of codimension one are: an Andronov–Hopf bifurcation with frequency of order $O(1)$ which we label $\mathrm{H}_1$, an Andronov–Hopf bifurcation with small frequency $O(|\alpha_1|)$, labeled $\mathrm{H}_0$, and a pitchfork bifurcation labeled P; see Figure 1. Recall that the unfolding parameters in Figure 1 are related to the physical parameters $\mu = (\Omega, \delta_1, \delta_2)$ by equations (2.12).

4.1.1. *Primary Andronov–Hopf Bifurcation* $\mathrm{H}_1$. For the system (3.3), the $z$–plane is a two-dimensional invariant subspace, on which the normal form equations reduce to

$$\begin{aligned} y_1 = y_2 &= 0 \\ \dot{r} &= \left[\alpha_3 + \phi_3(0, r^2)\right] r \\ \dot{\theta} &= 1 + \phi_4(0, r^2) \ . \end{aligned} \tag{4.2}$$

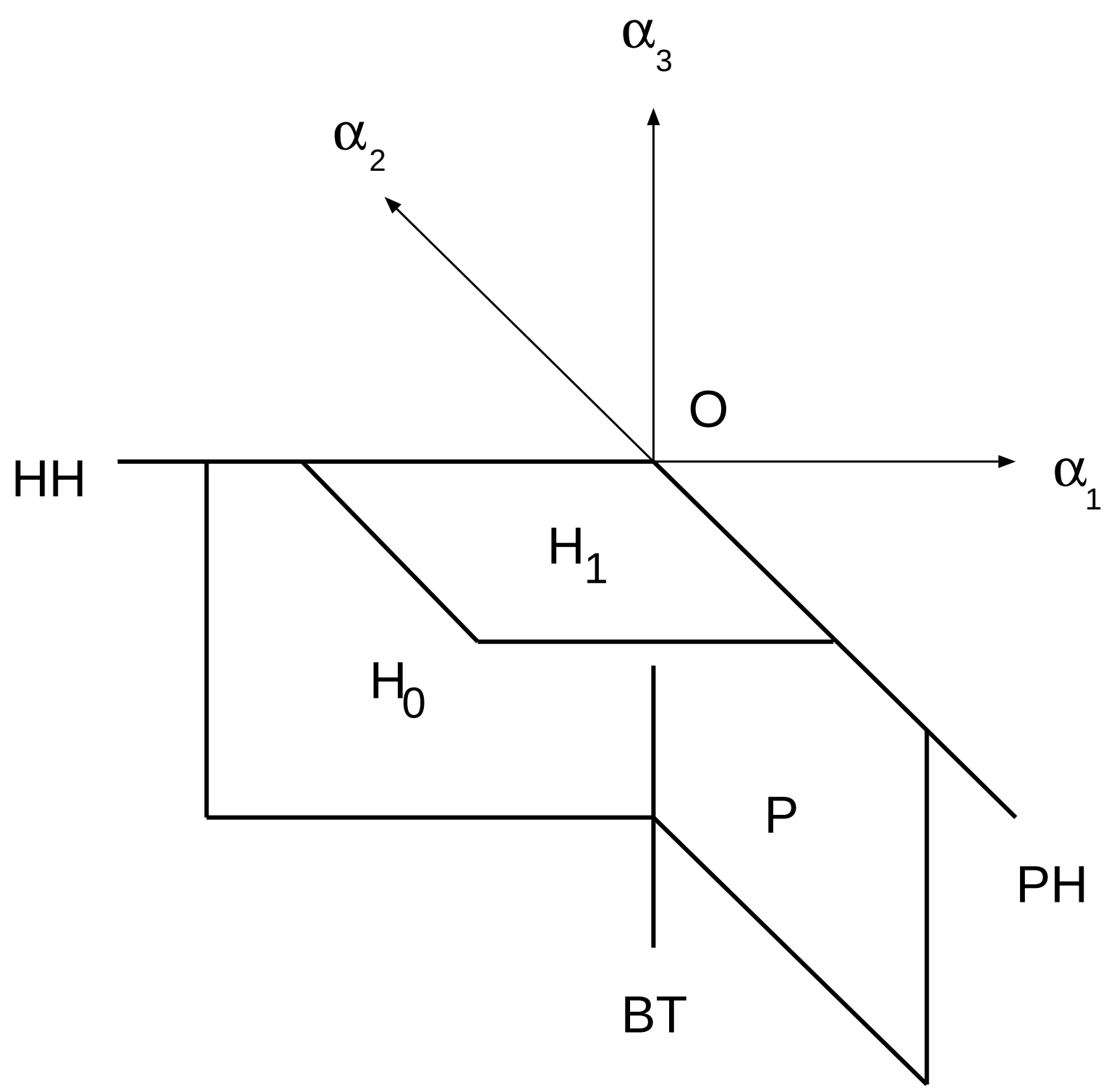

FIGURE 1. Stratified subvariety of primary bifurcations.

This has the family of periodic solutions

$$\begin{aligned} y_1 = y_2 &= 0 \\ r &= \sqrt{\tfrac{-\alpha_3}{\phi_{3,01}}} + \cdots \\ \theta(t) &= \theta_0 + \left[1 + \phi_4(0, r^2)\right] t \pmod{2\pi} . \end{aligned} \tag{4.3}$$

This is a classical Andronov–Hopf bifurcation, giving birth to a limit cycle in the $(y_3, y_4)$–plane with period near $2\pi$ (angular frequency 1) as the parameter $\alpha_3$ crosses through 0, provided $\alpha_1, \alpha_2 \neq 0$. Generically $\phi_{3,01} \neq 0$ and the bifurcation is nondegenerate, which is assumed to be the case in this paper; then the limit cycle is locally unique. Degenerate cases of Andronov–Hopf bifurcation are explored in detail in [**17**], but will not be considered here. If $\phi_{3,01} < 0$, then the bifurcation is *supercritical* (the limit cycle exists for $\alpha_3 > 0$) and if $\phi_{3,01} > 0$ the bifurcation is *subcritical* (exists for $\alpha_3 < 0$). A supercritical limit cycle is stable if $\alpha_1, \alpha_2 < 0$, otherwise the cycle is unstable.

4.1.2. *Pitchfork Bifurcation* P. Any equilibrium solution of (3.3) lies on the $y_1$ axis, and satisfies

$$\left[\alpha_1 + \phi_1(y_1^2, 0)\right] y_1 = 0 . \tag{4.4}$$

This is a classic *pitchfork bifurcation*, yielding in addition to the trivial solution $\mathbf{y} = 0$, a symmetric pair of nontrivial solutions

$$y_1 = \pm\sqrt{-\frac{\alpha_1}{\phi_{1,10}}} + \cdots$$

branching from 0 as $\alpha_1$ crosses through zero, provided $\phi_{1,10} \neq 0$ and $\alpha_2, \alpha_3 \neq 0$, which we assume. The pitchfork is supercritical if $\phi_{1,10} < 0$ and subcritical if $\phi_{1,10} > 0$. Supercritical solutions are stable if $\alpha_2, \alpha_3 < 0$, otherwise both are unstable.

4.1.3. *Primary Andronov–Hopf Bifurcation* $\mathrm{H}_0$. The Poincaré normal form (3.3) leaves the $(y_1, y_2)$–plane invariant. On this plane, the normal form reduces to

$$\begin{aligned} \dot{y}_1 &= y_2 \\ \dot{y}_2 &= [\alpha_1 + \phi_1(y_1^2, 0)]y_1 + [\alpha_2 + \phi_2(y_1^2, 0)]y_2 \\ z &= 0\,. \end{aligned} \tag{4.5}$$

For fixed $\alpha_1 < 0$, this system has an Andronov–Hopf bifurcation as $\alpha_2$ crosses through zero, for any $\alpha_3 \neq 0$. The new limit cycle has a small frequency, asymptotic to $\sqrt{-\alpha_1}$. The change of variables $x = y_1,\ y = -y_2/\omega$, where $\omega^2 = -\alpha_1$ transforms (4.5) to the planar system

$$\begin{pmatrix} \dot{x} \\ \dot{y} \end{pmatrix} = \begin{pmatrix} 0 & -\omega \\ \omega & 0 \end{pmatrix} \begin{pmatrix} x \\ y \end{pmatrix} + \begin{pmatrix} 0 \\ -\frac{1}{\omega}\phi_1(x^2, 0)x + \phi_2(x^2, 0)y \end{pmatrix}. \tag{4.6}$$

This planar system has the polar normal form, with $x + iy = \rho e^{i\varphi}$

$$\begin{aligned} \dot{\rho} &= \rho\left[\frac{\alpha_2}{2} + \frac{1}{8}\phi_{2,10}\rho^2 + \cdots\right] \\ \dot{\varphi} &= \sqrt{-\alpha_1 - \left(\frac{\alpha_2}{2}\right)^2} + b\rho^2 + \cdots\,. \end{aligned} \tag{4.7}$$

From this it follows that the limit cycle has amplitude given by

$$\rho = \sqrt{-\frac{4\alpha_2}{\phi_{2,10}}} + \cdots$$

and is super-(sub-)critical if $\phi_{2,10} < 0$ ($> 0$).

**4.2. Primary Bifurcations of Codimension Two.** The Poincaré normal form (3.3) has three bifurcations of codimension two, from the trivial solution in state space. In the parameter space these occur along the three coordinate axes, see Figure 1.

4.2.1. *Bogdanov–Takens Bifurcation* BT. As observed for equation (4.5), the $(y_1, y_2)$–plane is invariant for the flow of the normal form. On this plane, for $\alpha_3 \neq 0$ and $(\alpha_1, \alpha_2)$ in a neighbourhood of $(0, 0)$, the system has a $\mathbb{Z}_2$–symmetric Bogdanov–Takens bifurcation [**11, 41**], see line BT in Figure 1. It yields stable solutions only if $\alpha_3 < 0$. In addition to the codimension one primary bifurcations P and $\mathrm{H}_0$ described above, there are rich possibilities including simultaneous secondary Andronov–Hopf bifurcations from the pitchfork pair of equilibria, coalescence of two limit cycles, orbits homoclinic to the origin and heteroclinic orbits joining the two pitchfork equilibria. A detailed description of the versal deformations, bifurcation diagrams and phase portraits is given by Arnold [**9**, page 300], and need not be repeated here.

New possibilities arise in this system via the interaction of the Bogdanov–Takens bifurcation with the primary Andronov–Hopf bifurcation $\mathrm{H}_1$, leading to a flow in a phase space which reduces to $\mathbb{R}^2 \times \mathrm{S}^1$. The simple homoclinic orbits of the planar Bogdanov–Takens bifurcation then lead to transversally intersecting two-dimensional stable and unstable manifolds of periodic orbits and thus to chaotic dynamics, as for example in [**18, 28, 42, 44**]. These phenomena are under further investigation.

4.2.2. *Pitchfork–Hopf Mode Interaction* PH. On the $\alpha_2$–axis $(\alpha_1, \alpha_3) = (0, 0)$, the system has a simple zero eigenvalue, with a $\mathbb{Z}_2$–symmetry on the corresponding eigenspace $y_1$, and a pair of purely imaginary eigenvalues $\pm i$. In a neighbourhood of the $\alpha_2$–axis for $\alpha_2 < 0$, the normal form (3.3) has a pitchfork-Hopf mode interaction, see PH in Figure 1. One finds the expected pitchfork and Andronov–Hopf primary bifurcations P and $\mathrm{H}_1$, and secondary bifurcations of mixed-mode solutions involving nonlinear coupling of these primary bifurcations. A complete listing of these primary and secondary bifurcation branches for nondegenerate pitchfork–Hopf mode interactions is given in [**27**], see also [**18**]. There is a possibility also of tertiary bifurcation to a quasiperiodic flow on an invariant torus with a small second frequency [**25, 39**], discussed further below.

To investigate these mode-interactions, in (3.4) we set $\dot{y}_1 = y_2 = 0$, $\dot{r} = 0$, suppress the $\dot{\theta}$ equation and retain only leading order terms in $\phi_1$ and $\phi_3$. Then equation (3.5) reduces to

$$\begin{array}{rcl} 0 & = & \left[\alpha_1 + \phi_{1,10} y_1^2 + \phi_{1,01} r^2\right] y_1 \\ 0 & = & \left[\alpha_3 + \phi_{3,10} y_1^2 + \phi_{3,01} r^2\right] r \ . \end{array} \tag{4.8}$$

For nondegeneracy, we assume

$$\Delta_1 \equiv \phi_{1,10}\phi_{3,01} - \phi_{3,10}\phi_{1,01} \neq 0, \quad \phi_{1,10}\phi_{3,01} \neq 0 \ . \tag{4.9}$$

Then (4.8) has four branches of solutions

$$\begin{array}{rcl} (y_1^2, r^2) & = & (0, 0) \\ (y_1^2, r^2) & = & \left(-\dfrac{\alpha_1}{\phi_{1,10}}, 0\right) \\ (y_1^2, r^2) & = & \left(0, -\dfrac{\alpha_3}{\phi_{3,01}}\right) \\ (y_1^2, r^2) & = & \left(\dfrac{\alpha_3\phi_{1,01} - \alpha_1\phi_{3,01}}{\Delta_1}, \dfrac{\alpha_1\phi_{3,10} - \alpha_3\phi_{1,10}}{\Delta_1}\right) \ . \end{array} \tag{4.10}$$

Note that the nonzero quantities on the right must be positive for the solutions to be real. The first three of these solutions are the previously identified trivial, primary pitchfork P and primary Hopf $\mathrm{H}_1$ solutions, respectively. The fourth is a mixed-mode solution, which joins the P and $\mathrm{H}_1$ solutions, at the two ends of its interval of existence, in secondary bifurcations. Due to symmetry, the primary bifurcations have two symmetric branches $\pm y_1$, $\pm r$, and the secondary bifurcations have four branches; however, in the full dynamics the solutions with $\pm r$ give the same periodic orbit with a phase shift. The secondary bifurcation from P is an Andronov–Hopf bifurcation, while the secondary bifurcation from $\mathrm{H}_1$ is a pitchfork bifurcation of two new limit cycles. While all of these results were deduced only from the cubic truncation (4.8), they extend by standard arguments to the full equations, using assumptions (4.9), the Implicit Function Theorem and transversality.

When higher order terms are included, in the case $\phi_{1,10}\phi_{3,01} < 0$ it is possible to find a tertiary bifurcation, from the mixed-mode secondary branch above to a flow on a two-dimensional invariant torus with a small second frequency [**13, 18, 25, 39, 42**]. Explicit conditions for the existence of this type of 2-torus bifurcating from a steady state – Hopf mode interaction, and a description of the flow on the torus, have been given by Scheurle and Marsden, see [**39**]. It should be noted that this 2-torus exists only very locally in parameter space, in contrast to the 2-torus of the next subsection. The secondary radius of this 2-torus grows very rapidly from bifurcation, and it becomes sensitive to the effects of the symmetry–breaking tail (or remainder) of the normal form, which may lead to chaos.

4.2.3. *Hopf–Hopf Mode Interaction* HH. In a neighbourhood of the negative $\alpha_1$ semi-axis (see HH in Figure 1), nonresonant Hopf-Hopf mode interactions occur, involving the $\mathrm{H}_1$ and $\mathrm{H}_0$ bifurcations. For $(\alpha_2, \alpha_3) = (0, 0)$, the full system has two pairs of purely imaginary eigenvalues, which are nonresonant because the frequency ratio is very small. An introduction to the dynamics of nonresonant Hopf-Hopf mode interactions is given in [**18**, Chapter 7]. The four-dimensional system can be written as two pairs of amplitude-phase equations, in $(r, \theta)$ and $(\rho, \varphi)$, corresponding to the two primary Andronov–Hopf bifurcations $\mathrm{H}_1$ and $\mathrm{H}_0$ respectively. The two phase equations decouple due to nonresonance, leaving two amplitude equations which have the $\mathbb{Z}_2 \oplus \mathbb{Z}_2$ symmetric normal form

$$\begin{aligned} \dot{r} &= r\left[\alpha_3 + c\rho^2 + \phi_{3,01} r^2 + \cdots\right] \\ \dot{\rho} &= \rho\left[\tfrac{1}{2}\alpha_2 + \tfrac{1}{8}\phi_{2,10}\rho^2 + dr^2 + \cdots\right] \ . \end{aligned} \tag{4.11}$$

This has four families of solutions satisfying $(\dot{r}, \dot{\rho}) = (0, 0)$, given, subject to the obvious nondegeneracy conditions, to leading order by

$$\begin{aligned} (r,\rho) &= (0,0) \\ (r,\rho) &= \left(\sqrt{-\frac{\alpha_3}{\phi_{3,01}}}, 0\right) \\ (r,\rho) &= \left(0, \sqrt{-\frac{4\alpha_2}{\phi_{2,10}}}\right) \\ (r,\rho) &= \left(\sqrt{\frac{\frac{1}{2}\alpha_2 c - \frac{1}{8}\alpha_3\phi_{2,10}}{\Delta_2}}, \sqrt{\frac{\alpha_3 d - \frac{1}{2}\alpha_2\phi_{3,01}}{\Delta_2}}\right) \end{aligned} \tag{4.12}$$

where $\Delta_2 \equiv \frac{1}{8}\phi_{2,10}\phi_{3,01} - cd$. The first three of these are the previously obtained trivial, Andronov–Hopf $\mathrm{H}_1$ and Andronov–Hopf $\mathrm{H}_0$ solutions, respectively. The fourth gives a new family of two-frequency solutions, on an invariant 2-torus. The local persistence of this invariant torus is guaranteed by standard arguments involving normal hyperbolicity and nonresonance. This 2-torus differs from the 2-tori discussed in the previous subsection, in that it is more robust, and exists over a much wider region of parameter space. Thus, one may expect that this torus may be observed more easily in physical systems.

This 2-torus appears through a secondary bifurcation from either of the two primary Andronov–Hopf families, and may smoothly link these two primary families. Generically, this family of 2-tori is either always stable or always unstable. To leading order, the 2-torus bifurcates from the $H_1$ periodic orbit when

$$\alpha_3 d = \frac{1}{2}\alpha_2\phi_{3,01}, \ \rho = 0, \ r = \sqrt{-\frac{\alpha_3}{\phi_{3,01}}}, \tag{4.13}$$

and bifurcates from the $H_0$ periodic orbit when

$$\alpha_2 c = \frac{1}{4}\alpha_3\phi_{2,10}, \quad \rho = \sqrt{-\frac{4\alpha_2}{\phi_{2,10}}}, \ r = 0 \ . \tag{4.14}$$

For parameters near (4.13), the solutions on the 2-torus are close to the $H_1$ periodic orbit with primary frequency $\mathcal{O}(1)$, but have a small amplitude modulation with frequency $\mathcal{O}(\sqrt{-\alpha_1})$. The amplitude of the modulation grows along the branch from (4.13) towards (4.14). This behaviour resembles that observed by Nayfeh et al. [**34, 35, 36**]. In certain cases, there is a further possibility of a stable 3-torus bifurcating locally from the 2-torus [**21**], and of chaotic dynamics; however, the analysis of these phenomena requires consideration of higher order terms and will not be pursued here.

## 5. Example

The analysis of the previous section makes it easy to find an example of a pair of coupled oscillators, of the general form (1.3), which exhibits (for a range of parameter values) the behaviour described by Nayfeh et al. [**34, 35, 36**]. One such example is

$$\begin{aligned} \ddot{x} + \delta_1\dot{x} + \Omega x + (x + \dot{x})(x^2 + y^2) &= 0 \\ \ddot{y} + \delta_2\dot{y} + y + (y + \dot{y})(-0.2x^2 + y^2) &= 0 \end{aligned} \tag{5.1}$$

with parameter values chosen as $\Omega = 0.3$, $\delta_1 = -0.2$ and $\delta_2 = -0.25$. Equations (5.1) may be interpreted as a pair of nonlinearly coupled Duffing – van der Pol oscillators, as follows. In the first equation, the $x^3$ term is the Duffing nonlinear restoring force and $\dot{x}x^2$ is the van der Pol nonlinear damping force. These are both modified through nonlinear coupling effects of the second oscillator, represented by the $xy^2$ and $\dot{x}y^2$ terms, respectively. An analogous interpretation holds for the second equation.

In accord with the analysis of the previous section, solutions of this system approach an invariant 2-torus, which exhibits a “fast oscillation” in the $(y, \dot{y})$ variables together with a “slow modulation” in the $x$ variable; see Figure 2 for two such solutions, which have been computed numerically and plotted using Maple. The physical parameters chosen in this example correspond to values of the unfolding parameters near the negative $\alpha_1$-axis in Figure 1, given by equation (2.12) approximately as

$$\alpha_1 = -0.3, \ \alpha_2 = 0.2, \ \alpha_3 = 0.125\,. \tag{5.2}$$

A more detailed study of this phenomenon, for more physically relevant examples, and including both numerical and analytical results, is in progress.

## 6. Conclusions

We may draw several conclusions from this study. The first important conclusion is that the simultaneous occurrence of Andronov–Hopf and Bogdanov–Takens bifurcations is a codimension-three phenomenon, requiring three independent parameters for its exploration and understanding. In the versal deformation of the codimension three bifurcation, there are primary bifurcations of six different types; three of codimension one and three of codimension two, all of which are classical. The stratified subvariety of primary bifurcations is shown in Figure 1. The three parameters in Arnold's versal deformation may be chosen so that they directly determine each of the three primary bifurcations of codimension one: the plane $\alpha_1 = 0$ determines the primary pitchfork bifurcation P, $\alpha_2 = 0$ determines the primary Andronov–Hopf bifurcation $\mathrm{H}_0$ with low frequency (of order $\sqrt{-\alpha_1}$) and $\alpha_3 = 0$ determines the primary Andronov–Hopf bifurcation $\mathrm{H}_1$ with high frequency (of order 1). The three codimension two bifurcations correspond to the three coordinate axes in Figure 1. They lead to secondary bifurcations of new families of solutions which are linkages between the primary bifurcation families of codimension one, and to further rich and complex dynamics.

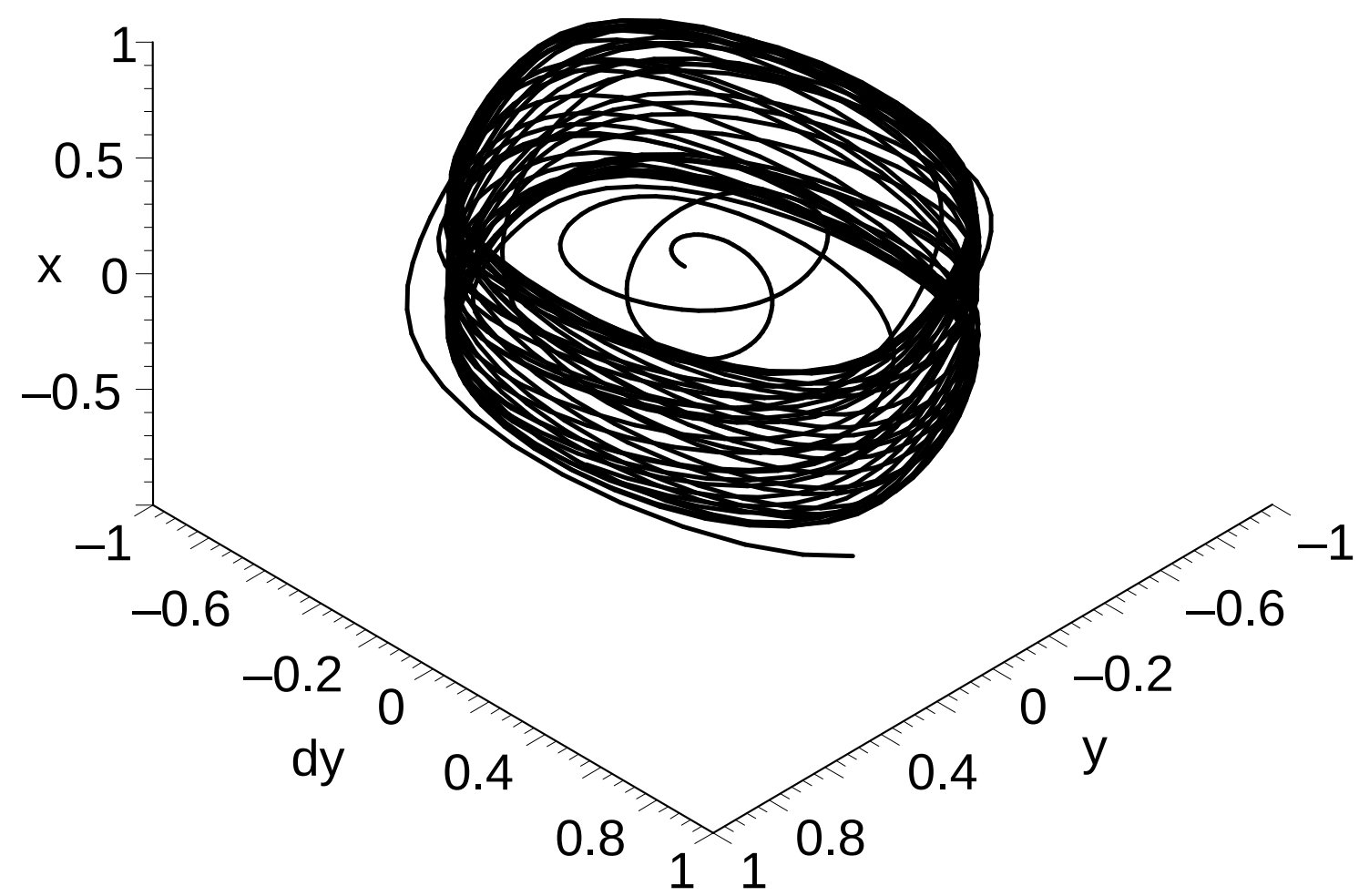

FIGURE 2. Two solutions approaching the invariant torus of Section 5, from two initial points inside and outside the torus, $(x(0), \dot{x}(0), y(0), \dot{y}(0)) = (0.1, 0, 0.1, 0)$ and $(-0.3, 0, 0.4, 0.9)$, respectively.

The three versal unfolding parameters are related by a nonsingular transformation to the three physical parameters of the coupled oscillators in the Introduction, by formulae (2.12). This verifies that a versal deformation gives the proper setting for investigation of the possible behaviours of the class of coupled oscillators defined in the Introduction.

The nonlinear coupling between the primary high frequency and low frequency Andronov–Hopf bifurcations can lead to a quasiperiodic flow on a 2-torus, with its two frequencies close to those of the two primary Andronov–Hopf bifurcations. This gives a plausible explanation for the observations of Nayfeh et al. [**34, 35, 36**], of low frequency modulation of a high frequency mode, and of energy transfer from high to low frequency modes, near 0 : 1 resonance. The case studied by Nayfeh corresponds to a choice of unfolding parameters $\alpha_1 < 0$ (fixed) and $\alpha_2, \alpha_3$ both small and negative; i.e., near the negative $\alpha_1$ axis marked HH in Figure 1. This is the region where existence of a 2-torus is predicted by the normal form analysis and confirmed by the numerical analysis in Section 5.

The present study yields several indicators of directions for future research. Systems of coupled oscillators with less symmetry than required by (1.2) are important in applications, and may be investigated by the same normal form methods as used here. The effect of higher order terms, particularly in the onset of more complex dynamics, is worthy of investigation. For example, in the standard symmetric Bogdanov–Takens analysis, the double zero eigenvalue may lead to homoclinic and heteroclinic orbits in the phase plane. When coupled with the primary Andronov–Hopf bifurcation $H_1$ as here, it will lead to phenomena of homoclinic chaos, which merit further investigation. A more detailed investigation of the mechanism of energy transfer from high to low frequency modes, for more realistic physical examples, and using both analytical and numerical techniques, is in progress. As well, the application of these results to concrete systems in engineering, biology and other fields may lead to new insights.

## References

[1] Andronov, A.A. [1929] Les cycles limits de Poincaré et la théorie des oscillations auto-entretenues. *C.R. Acad. Sci.* **189**, 559–561.

[2] Andronov, A.A., Leontovich, E.A., Gordon, I.I. and Maier, A.G. [1973] *Qualitative Theory of Second-Order Dynamic Systems.* (Translation of the Russian edition of 1966). Halsted Press, John Wiley & Sons, New York.

[3] Andronov, A.A., Leontovich, E.A., Gordon, I.I. and Maier, A.G. [1973] *Theory of Bifurcations of Dynamic Systems on a Plane.* (Translation of the Russian edition of 1967). Halsted Press, John Wiley & Sons, New York.

[4] Andronov, A.A., Vitt, A.A. and Khaikin, S.E. [1966] *Theory of Oscillators.* (Translation of the second Russian edition). Pergamon Press Ltd., Oxford.

[5] Arnold, V.I. [1965] Small denominators I. Mappings of the circumference onto itself. *Amer. Math. Soc. Transl. Ser. 2.* **46**, 213–284.

[6] Arnold, V.I. [1971] On matrices depending on parameters. *Russian Math. Surveys* **26**, 29–43.

[7] Arnold, V.I. [1972] Lectures on bifurcations in versal families. *Russian Math. Surveys* **27**, 54–123.

[8] Arnold, V.I. [1977] Loss of stability of self-oscillations close to resonance and versal deformations of equivariant vector fields. *Func. Anal. Applic.* **11**, 85–92.

[9] Arnold, V.I. [1983] *Geometrical Methods in the Theory of Ordinary Differential Equations.* Springer-Verlag, New York.

[10] Bibikov, Y.N. [1979] *Local Theory of Nonlinear Analytic Ordinary Differential Equations.* LNM 702, Springer-Verlag, New York.

[11] Bogdanov R.I. [1975] Versal deformations of a singular point on the plane in the case of zero eigenvalues. *Funct. Anal. Appl.* **9**, 144–145.

[12] Bruno, A.D. [1989] *Local Methods in Nonlinear Differential Equations.* Springer-Verlag, Berlin.

[13] Broer, H.W. [1981] Quasiperiodic flow near a codimension one singularity of a divergence free vector field in dimension three. *Lect. Notes Math.* **898**, Springer-Verlag, New York.

[14] Chow, S.-N., Li, C. and Wang, D. [1994] *Normal forms and bifurcations of planar vector fields.* Cambridge University Press.

[15] Georgiou, I.T., Bajaj, A.K. and Corless, M. [1998] Slow and fast invariant manifolds, and normal modes in a two-degree-of-freedom structural dynamical system with multiple equilibrium states. *Internat. J. Non-Linear Mech.* **33**, 275–300.

[16] Gils, S.A. van, Krupa, M.P. and Langford, W.F. [1990] Hopf bifurcation with nonsemisimple 1:1 resonance. *Nonlinearity* **3**, 1–26.

[17] Golubitsky, M. and Langford, W.F. [1981] Classification and unfoldings of degenerate Hopf bifurcations. *J. Diff. Eq.* **41**, 375–415.

[18] Guckenheimer, J. and Holmes, P. [1983] *Nonlinear Oscillations, Dynamical Systems and Bifurcations of Vector Fields. Appl. Math Sci.* **42**, Springer-Verlag, New York.

[19] Holmes, P., Marsden, J.E. and Scheurle, J. [1988] Exponentially small splittings of separatrices with applications to KAM theory and degenerate bifurcations. *Hamiltonian Dynamical Systems (Boulder CO 1987) Contemp. Math.* **81** 213–244, Amer. Math. Soc., Providence RI.

[20] Iooss, G. and Adelmeyer, M. [1992] *Topics in Bifurcation Theory and Applications.* World Sci., Singapore.

[21] Iooss, G. and Langford, W.F. [1980] Conjectures on the routes to turbulence via bifurcations. In *Nonlinear Dynamics, Annals N.Y. Acad. Sci.* **357**, 489–505.

[22] Kahn, P.B. and Zarmi, Y. [1998] *Nonlinear Dynamics: Exploration through Normal Forms.* Wiley-Interscience, New York.

[23] Khazin, L.G. and Shnol, E.E. [1991] *Stability of Critical Equilibrium States.* Manchester Univ. Press, Manchester.

[24] Langford, W.F. [1979] Periodic and steady mode interactions lead to tori. *SIAM J. Appl. Math.* **37**, 22–48.

[25] Langford, W.F. [1983] A review of interactions of Hopf and steady-state bifurcations. In *Dynamical Systems and Turbulence.* G.I. Barenblatt, G. Iooss and D.D. Joseph editors, Pitman, London, 215–237.

[26] Langford, W.F. [1984] Numerical studies of torus bifurcations. In *Numerical Methods for Bifurcation Problems.* T. Küpper, H.D. Mittelmann and H. Weber editors, ISNM 70 Birkhäuser Verlag, Boston, 285–294.

[27] Langford, W.F. and Iooss, G. [1980] Interactions of Hopf and pitchfork bifurcations. In *Bifurcation Problems and their Numerical Solution.* H.D. Mittelmann and H. Weber editors, ISNM 54, Birkhäuser Verlag, Basel, 103–134.

[28] Langford, W.F. and Nagata, W. [1995] *Normal Forms and Homoclinic Chaos, Fields Institute Communications* **4**, American Mathematical Society, Providence.

[29] Luo, D., Wang, X., Zhu, D. and Han, M. [1997] *Bifurcation Theory and Methods of Dynamical Systems.* World Scientific, Singapore.

[30] Marsden, J.E. and McCracken, M. [1976] *The Hopf Bifurcation and its Applications, Appl. Math. Sci.* **19**, Springer–Verlag, New York.

[31] Moson, P. [1984] Quasi-periodic solutions of a special system. In *Differential Equations: Qualitative Theory.* North-Holland, Amsterdam.

[32] Moson, P. [1991] Local bifurcations in the case of eigenvalues $0^2$, $\pm i$. *Z. Angew. Math. Mech.* **71**, T69–T70.

[33] Nagata, W. and Namachchivaya, N. Sri [1998] Bifurcations in gyroscopic systems with an application to rotating shafts. *Proc. Roy. Soc. London Ser. A* **454**, 543–585.

[34] Nayfeh, S.A. and Nayfeh, A.H. [1993] Nonlinear interactions between two widely spaced modes – external excitation. *Internat. J. Bifur. Chaos* **3**, 417–427.

[35] Nayfeh, S.A. and Nayfeh, A.H. [1994] Energy Transfer from high- to low-frequency modes in a flexible structure via modulation. *J. Vib. Acoustics* **116**, 203–207.

[36] Nayfeh, A.H., Nayfeh, S.A., Anderson, T.A. and Balachandran, B. [1994] Transfer of energy from high-frequency to low-frequency modes. In *Nonlinearity and Chaos in Engineering Dynamics.* John-Wiley & Sons, Chichester.

[37] Poincaré, H. [1993] *New Methods of Celestial Mechanics.* (Translation of *Les Méthodes nouvelles de la Méchanique céleste*, originally published in 1892–1899). Three volumes, edited and introduced by D.L. Goroff, American Institute of Physics, College Park, MD.

[38] Sanders, J.A. [1982] Melnikov's method and averaging. *Celestial Mechanics* **28** 171–181.

[39] Scheurle, J. and Marsden, J.E. [1984] Bifurcation to quasi-periodic tori in the interaction of steady state and Hopf bifurcations. *SIAM J. Math. Anal.* **15** 1055–1074.

[40] Scheurle, J., Marsden, J.E. and Holmes, P. [1991] Exponentially small estimates for separatrix splittings. *Asymptotics Beyond all Orders (La Jolla CA 1991) NATO Adv. Sci. Inst. Ser. B Phys.* **284** 187–195, Plenum, New York.

[41] Takens, F. [1974] Singularities of vector fields. *Publ. Math. IHES* **43**, 47–100.

[42] Wiggins, S. [1988] *Global Bifurcations and Chaos: Analytical Methods. Appl. Math Sci.* **73**, Springer-Verlag, New York.

[43] Wiggins, S. [1990] *Introduction to Applied Nonlinear Dynamical Systems and Chaos. Texts Appl. Math.* **2**, Springer-Verlag, New York.

[44] Wiggins, S. [1993] *Global Dynamics, Phase Space Transport, Orbits Homoclinic to Resonances, and Applications. Fields Inst. Monographs* **1**, Amer. Math. Soc., Providence RI.
[45] Yu, P. and Huseyin, K. [1988] Bifurcations associated with a double zero and a pair of pure imaginary eigenvalues. *SIAM J. Appl. Math.* **48**, 229–261.

Department of Mathematics and Statistics, University of Guelph, Guelph, Ontario, Canada N1G 2W1
*E-mail address*: `wlangfor@uoguelph.ca`